\documentstyle{amsppt}
\magnification =\magstep 1
\hsize 6.5 true in
\vsize 8.9 true in
\NoBlackBoxes
\topmatter
\title Operations on Cyclic Homology, the $X$ Complex, and a
Conjecture of Deligne \endtitle
\abstract We prove that there is a product on the Hochschild and cyclic
 chain complex of a homotopy Gerstenhaber algebra. By restricting to
  the special case of
 algebra of Hochschild cochains (the so called deformation complex), we
 obtain operations on cyclic homology of associative algebras.
 \endabstract
\author Masoud Khalkhali \endauthor
\TagsOnLeft
\address Department of Mathematics \newline
University of Western Ontario \newline
London ON Canada \newline
N6A 5B7
\endaddress
\thanks{Supported in part by NSERC of Canada. I would also like to thank
the mathematics division of the
International Centre for Theoretical Physics, Trieste, Italy,  for a visiting fellowship
during the summer of 1996 where part of this work was completed.}
\endthanks
\endtopmatter
\rightheadtext{Operations on cyclic homology ...}
\document
\subhead  1\quad Introduction \endsubhead

The goal of this article is to relate  recent developments in cyclic
homology theory [3] and the theory of operads and homotopical algebra
[6,8], and hence
to provide a general framework to define and study operations in cyclic 
homology theory.
The link here is the bar construction.

In [4], P. Deligne conjectured that the Hochschild cochain complex of an
associative algebra, also called the deformation complex of the algebra, has a natural structure of an algebra over a
singular chain operad of the little squares operad. This conjecture is
now
proved by M. Kontsevich [16]. It should be noted, however, that
the results of the present paper in no way depends on the geometric form
of this conjecture. In
fact we do not use the conjecture in its original geometric form.   A closely related statement recently proved by
Gerstenhaber and Voronov [6] and Getzler
and Jones [8] states that the deformation complex of an
associative algebra has a natural structure of a homotopy Gerstenhaber
algebra, also called homotopy $G$
algebra. This result completes Gerstenhaber's  earlier work in [5]
in the sense that it
reveals the full structure of higher homotopies in the deformation
complex. In a sense, this result shows that there is a natural "quantum group"
structure on the bar construction of the deformation complex. It is this algebraic version of the conjecture that is most
useful to construct the operations.

In fact, we found it  more conceptual to go beyond the deformation complex for
associative
algebras and define certain operations on the Hochschild and cyclic
complexes of
homotopy $G$ algebras (Theorem 8 below). As an application,  by
specializing to the homotopy $G$ algebra structure of the deformation
complex, we obtain  operations similar to those constructed by  Nest and
Tsygan in their study of
algebraic index theorems [12]. In view of increasing importance of homotopy
$G$ algebras and the theory of operads in general in mathematics and mathematical physics
(see, for example,  [6,8,11,15] and references therein),
we hope that theorem 8, and specially its  method of proof, which is
non-computational and lends itself to generalizations to algebras over
operads,
 will prove useful in applcations of noncommutative
geometry and cyclic homology.

This paper is organized as follows. In Sect. 2 we recall the notion of
homotopy $G$ algebra and specially its formulation in terms of
the bar construction from [6]. In Sect. 3  we
define operations on Hochschild and cyclic complex of homotopy $G$ algebras.
A central tool here is the
notion of $X$ complex and its refinements due to Cuntz and Quillen [3].
By a result of Quillen [14], cyclic and Hochschild chain complexes appear
as the $X$ complex of the bar construction. From this point of view
operations on cyclic and Hochschild complex of homotopy $G$ algebras
are predicted by  Kunneth formulas for
the $X$ complex of differential graded coalgebras. Sect. 4 is mainly devoted to
deriving explicit  formulas in the context of Connes' $b, B$ bicomplex
and also specializing to
the case of deformation complex.

I am much obliged to Maxim Kontsevich for a very informative
communication on Deligne's conjecture.

\subhead  2\quad Homotopy Gerstenhaber algebras \endsubhead

This section is based on [6]. In an attempt to make the paper as self
contained as possible, we have reproduced the proofs of the main
statements. Let $V=\bigoplus_{i\in \Bbb Z} V_i$ be a $\Bbb Z$-graded linear space. We use $\vert
x\vert $ to denote the degree of a homogeneous element $x\in V$. A {\it
brace algebra structure}  on $V$ is given by a collection of linear
homogeneous maps of
degree zero, indexed by $n\geq 0$,
$$V\otimes V^{\otimes n}\longrightarrow V$$
$$ x\otimes x_1\otimes \cdots \otimes x_n \mapsto x\{x_1,\cdots , x_n 
\},$$
such that, for all $m$,$n$,  the following {\it higher pre-Jacobi identities}  
are
satisfied:
$$ x\{x_1,\cdots ,x_m\}\{y_1,\cdots , y_n\}=
\sum_{0\leq i_1\leq j_1\leq \cdots \leq i_m\leq j_m\leq n}
(-1)^{\epsilon} \tag 1$$
$$x
\{y_1,\cdots ,y_{i_1},x_1\{y_{i_1+1},\cdots 
,y_{j_1}\},\cdots ,x_m\{y_{i_m+1},\cdots ,y_{j_m}\},\cdots,   y_n\}, $$
where $x_i$ and $y_j$ are homogeneous elements and $\epsilon =
\sum_{p=1}^{m} (\vert x_p\vert \sum_{q=1}^{i_p}\vert
y_q \vert).$ The degree zero assumption simply means that $\vert
x\{x_1,\cdots ,x_n\}\vert =\vert x\vert +\sum \vert x_i\vert $. We also
assume that for $n=0$ the resulting map $x\mapsto x\{ \} 
:V\longrightarrow V$ is the identity.

For example, as a consequence of (1),  one checks that the bracket
$$[x,y]:=x\{y\}-(-1)^{\vert x\vert \vert y\vert }y\{x\}$$
defines a graded
Lie algebra structure on $V$. Indeed, puting $m=n=1$ in (1), one obtains
$$x\{y\}\{z\}-x\{y\{z\}\}=x\{y,z\}+(-1)^{\vert y\vert \vert
z\vert}x\{z,y\},$$ 
which measures the failure of the operation $(x,y)\mapsto x\{y\}$ to be 
associative. From
this the graded Jacobi identity easily follows.

A brace algebra structure on $V$ has a particularly simple interpretation in
terms of the  {\it tensor coalgebra} of $V[1]$.   Let $V[1]$ be the desuspension of $V$
defined by $V[1]_n=V_{n+1}$. Let
$$T(V[1])=\bigoplus_{n\geq 0} (V[1])^{\otimes
n}$$
be the tensor coalgebra of $V[1]$ with its coproduct
$$\Delta
:T(V[1])\longrightarrow T(V[1]) \otimes T(V[1]) $$
$$\Delta (x_1,\cdots , x_n)=\sum_{i=0}^{n}(x_1,\cdots ,x_i)\otimes 
(x_{i+1}, \cdots ,x_n),$$
where
we have denoted a tensor $x_1\otimes
\cdots \otimes x_n$ in $T(V[1])$ by $(x_1, \cdots ,x_n)$.

Note that $T(V[1])$ is
 bigraded. Its horizontal grading
is denoted by $deg $ and
defined by $deg (x_1,\cdots, x_n)=n$, and its vertical grading, denoted
$\vert \quad \vert_d $, is
given by $\vert (x_1,\cdots ,x_n)\vert_d =\sum \vert x_i \vert -n$. The
total grading is hence given by $\vert (x_1,\cdots x_n)\vert =\sum_i
\vert x_i \vert$. In the definition of $T(V[1])$ we are implicitly assuming
that we are working with the total complex which is ${\Bbb Z}$-graded.

A linear homogeneous (with respect to the total grading) map
of total degree zero,
$$\cup :T(V[1])\otimes T(V[1]) \longrightarrow T(V[1]),$$
is called {\it left
increasing} if $deg (\alpha \cup \beta )\geq deg (\alpha)$. We always
assume that $\cup $ is counital.

\proclaim{Lemma 1} There is a natural $1-1$ correspondence between brace algebra
structures on $V$ and left increasing bialgebra structures on $T(V[1])$.
\endproclaim
\demo{Proof} Since $T(V[1])$ is the free coalgebra
generated by the desuspension $V[1]$, we have a natural $1-1$ correspondence
between coalgebra morphisms
                $$\cup :T(V[1])\otimes T(V[1]) \longrightarrow T(V[1])$$
                and linear maps
                $$m:T(V[1])\otimes T(V[1]) \longrightarrow V[1] .$$
Given $m$, the $n$-th component of $\cup $ is defined by
$$\cup_n=m^{\otimes n}\circ \tilde{\Delta}^{(n-1)}, \tag 2$$
where $\tilde{\Delta}^{(n)}$ denotes  the $n-th$ iteration of the
coproduct $\tilde{\Delta}$
of the coalgebra $T(V[1])\otimes T(V[1])$. Conversely, given $\cup $, $m$ is just the degree one
component of $\cup$.

It is clear that $\cup$ is left increasing iff, for $n\geq 2$, $m\vert
V[1]^{\otimes n}\otimes T(V[1])=0$. In this case, let us denote the map
$m : V[1]\otimes T(V[1])\longrightarrow V[1]$ by $x\{x_1,\cdots ,x_n\}.$
Then, using
(2), we get  an explicit formula for the coalgebra map $\cup :T(V[1])\otimes
T(V[1])\longrightarrow T(V[1])$. It is given by
$$(x_1,\cdots ,x_m)\cup (y_1,\cdots ,y_n)=\sum_{0\leq i_1\leq j_1\leq
\cdots \leq i_m\leq j_m \leq n} (-1)^{\epsilon} \tag 3 $$
$$(y_1,\cdots y_{i_1},x_1\{y_{i_1+1},\cdots ,y_{j_1}\},\cdots 
,x_m\{y_{i_m +1},\cdots ,y_{j_m}\},\cdots ,y_n),$$
where  $\epsilon $ is the same as in (1).

It remains to check that $\cup$ is associative iff the braces  
satisfy the higher
pre-Jacobi identities (1).  Assume $\cup$ is associative. Then in particular we
have $(\alpha \cup \beta)\cup \gamma =\alpha \cup (\beta \cup \gamma)$, 
where $\alpha
=x$, $\beta =(x_1,\cdots ,x_m)$ and $\gamma =(y_1,\cdots y_n)$.
Taking the  degree one component of both sides,  one obtains (1).

Conversely, assume the braces  satisfy the pre-Jacboi identities. It
is easy to see that both maps
$ \cup \otimes 1$ and $1\otimes \cup :
T(V[1])\otimes T(V[1])\otimes T(V[1]) \longrightarrow T(V[1])$
are coalgebra maps. By the universal property of $T(V[1])$, the
two maps are the same provided the degree one components of them 
coincide.One then checks that this is equivalent to the brace
identity (1). The theorem is proved.
\enddemo

In the rest of this paper we only consider left increasing 
multiplications on $T(V[1])$.
Here is an example of a brace algebra. This example is due to Getzler
[7].  Let $A$ be a linear
space and let $V_n=Hom (A^{\otimes n}, A)$. One defines a brace algebra
structure on $V[1]$ by setting
$$x\{x_1,\cdots ,x_m\}(a_1,\cdots a_n)= \sum_{0\leq i_1\leq \cdots
\leq i_m\leq n}(-1)^{\epsilon} \tag 4 $$
$$x(a_1, \cdots ,x_1(a_{i_1+1},\cdots, a_{i_1+d_1} ),
\cdots x_m(a_{i_m+1},\cdots ,a_{i_m+d_m}
),\cdots , a_n),$$
where $d_i=\vert x_i \vert +1$, $n=1+\vert x \vert +\sum \vert x_i \vert $ and $\epsilon =
\sum_{p=1}^{m}\vert x_p\vert i_p  .$
Checking (1) is
straightforward.

A {\it homotopy G algebra} (G stands for Gerstenhaber) is a
differential graded ($DG$) associative algebra  equipped with a
system of "higher homotopies" so that its
cohomology is a {\it graded poisson algebra}.
More precisely, Let $(V,\delta )$ be a $DG$ algebra where we assume
the differential has degree $+1$. A homotopy $G$ algebra structure on
$V$ is given by a brace algebra structure on the desuspension
$V[1]$ of $V$  such that  the
following axioms are satisfied:
$$ (x_1x_2)\{y_1,\cdots ,y_n\}=\sum_{k=0}^{n}(-1)^{\epsilon}x_1\{y_1,\cdots
,y_k\}x_2\{y_{k+1},\cdots ,y_n\},\tag 5$$
 where $\epsilon =(\vert x_2\vert -1)(\vert y_1\vert +\cdots +\vert y_k
 \vert
 -k)$, and
$$
\align
\delta (x\{x_1,\cdots ,x_{n+1}\})-\delta x\{x_1,\cdots ,x_{n+1}\} &-
\tag 6\\
(-1)^{\vert x\vert -1}\sum_{i=1}^{n+1}
(-1)^{\epsilon}x\{x_1,\cdots ,\delta x_i,\cdots ,x_{n+1}\} &=\\
-(-1)^{(\vert x\vert-1)(\vert x_1 \vert -1)}
x_1 .x\{x_2,\cdots ,x_{n+1}\}&+ \\
(-1)^{\vert x \vert -1}
\sum_{i=1}^n (-1)^{i+n}x\{x_1,\cdots
,x_ix_{i+1},\cdots ,x_{n+1}\} & \\
-x\{x_1,\cdots ,x_n\}.x_{n+1} &
\endalign
$$
 where $\epsilon =\sum_{k=1}^{i-1}\vert x_k\vert $.

A large class of homotopy $G$  algebras are constructed as follows [6].
Let $V$ be a brace algebra and $m\in V_1$ a degree one element
such that $m\{m\}=0$. One defines a $DG$ algebra  structure on $V[-1]$ by
defining a differential and a product by
$$\delta x=(-1)^{\vert  x\vert}[m,x]\quad \quad \quad \quad \quad
xy=m\{x,y\} \tag 7$$
Using the brace relations (1), one then
checks  that the axioms of a homotopy $G$ algebra are
satisfied.

The axioms of a  homotopy $G$ algebra structure on $V$ can be
conceptually encoded in terms of the {\it bar construction} $BV$. Let us 
describe this correspondence. We first need a definition.

Let $(V,\delta )$ be a $DG$ algebra where we assume the differential has
degree $+1$. Recal that the bar construction of $V$, denoted $BV$, is a
 differential graded coalgebra whose underlying coalgebra is the tensor
 coalgebra $T(V[2])$ and
its differential is the total differential $b'+\delta$. The individual
  differentials $ b',\delta  :BV \longrightarrow BV$
 are defined by
 $$b' (x_1,\cdots
 ,x_n)=(-1)^n\sum_{i=1}^{n-1}(-1)^{i-1}(x_1,\cdots,x_ix_{i+1},\cdots
 ,x_n),$$
$$ \delta (x_1,\cdots ,x_n)=\sum_{i=1}^n (-1)^{\vert x_1\vert +\cdots
+\vert x_{i-1}\vert } (x_1,\cdots ,\delta x_i ,\cdots
,x_n).$$
Note that both $b'$ and $\delta $ have total degree $+1$. Also note that
the total degree of $\alpha =(x_1,\cdots ,x_n) \in BV$ is given by
$\vert \alpha \vert =\sum_i \vert x_i \vert -n$.

Next recall
that a  $DG$ bialgebra is by definition a bialgebra object in the
abelian tensor category of
$DG$ linear spaces (cochain complexes). In particular the  differential
of a $DG$-bialgebra is simultaneously a derivation and a
coderivation.

\proclaim{Lemma 2} Let $V$ be a $DG$ algebra. Then there is a natural 1-1
correspondence between homotopy $G$ algebra structures on $V$ and $DG$
bialgebra structures on the bar construction $BV$.
\endproclaim

\demo{Proof} By the above lemma, we have a natural  $1-1$ correspondence 
between brace algebra structures on $V[1]$ and bialgebra structures on
$BV$. So
all that we need to prove is
that the axioms of  homotopy $G$ algebras (5 , 6) are
equivalent to  differential $b'+\delta $ being a derivation of $BV$.
That is for all $\alpha ,\beta \in BV$,
     $$(b'+\delta )(\alpha \cup \beta )=(b'+\delta)\alpha \cup \beta
     +(-1)^{\vert \alpha \vert }\alpha \cup
     (b'+\delta )\beta \tag 8$$
 
Now, since both $b'$ and $d$ are coderivations of $BV$ and $\cup $
is a coalgebra map, it follows that (8) holds if and only if the degree
one components of both sides coincide. 
Note that the only possible contributions to degree one components are 
from the following two choices:
$\alpha =x $, $\beta =(x_1,\cdots ,x_n)$ and $\alpha
=(x_1, x_2) $, $\beta =(y_1,\cdots ,y_m)$. Computing the first order
terms in the expansions, we find that (8) is equivalent to
(5, 6). The lemma is proved.
\enddemo

Given a homotopy $G$ algebra $V$, let $H(V)$ denote the cohomology of
the complex
$(V, \delta)$. The  product and the Lie bracket in $V$, being compatible with
the differential $\delta$, descend to $H(V)$ and define an associative
product and a Lie
algebra structure on $H(V)$. Moreover, the homotopy formulas in (5) and
(6) can
be used to show that the associative product in $V$ is, up to homotopy, 
graded commutative and
the Lie  bracket  is, again up to homotopy,  a derivation with
respect to the assocoative product. It thus follows that  $H(V)$
is a   graded poisson
algebra,  also known as a {\it Gerstenhaber algebra} ($G$ algebra). This
simply means that the product in cohomology is graded commutative and 
the Lie bracket is a derivation with respect to the product.

Examples of graded poisson algebras and homotopy $G$ algebras abound in
algebraic topology, Geometry and mathematical physics. By a classical
result of F. Cohen the cohomology groups of configuration spaces is a
universal model for graded poisson algebras in the sense that any graded
poisson algebra is an algebra over the latter as an operad. concrete
examples of $G$ algebras include the semi-infinte cohomology of string
theory [11], the algebra of polyvector fields on a
manifold, Koszul complex of Lie algebras and more generaly the
deformation cohomology of any associative algebra, to be discussed in
more detail in the next paragraph. Examples of homotopy $G$ algebra
structures that are just emerging include the homotopy $G$ algebra
structure of topological field theory [11], the homotopy $G$ algebra structure on
singular cochains on a topological space [6] and finally the  deformation
complex of associative algebras which we describe next. This structure
also appears in the recent work of M. Kontsevich on deformation
quantization of Poisson manifolds [15].

Let $A$ be an associative algebra and let $C (A,A)$ denote the
{\it deformation complex} of $A$. This is the standard
complex that calculates the Hochschild cohomology $H^{\bullet}(A,A)$.
We have $C^n(A,A)=Hom (A^{\otimes n},A)$. It thus follows from (4) that there
is a brace algebra structure defined on $C (A,A)[1]$.
Let $m:A\otimes A\longrightarrow A$ be the multiplication of $A$.
One has $m\{ m\} =0$, which is equivalent to associativity of
$m$. It is easy to check that the differential and the product induced
on $C(A,A)$
by (7) coincide, respectively, with the Hochschild coboundary and the
cup product on $C (A,A)$. One thus obtains a homotopy $G$ algebra
structure on $C (A,A)$, first discovered by Gerstenhaber and Voronov
in [6]. As we will see in the next section, this homotopy $G$ algebra
structure is at the heart of operations on Hochschild and cyclic
homology.

\subhead 3\quad Operations on homotopy $G$ algebras \endsubhead

Let $C$ and $D$ be $DG$ coalgebras and let $A$ be an algebra.
Our goal in this section is to show that any morphism of $DG$ coalgebras
$C\otimes D \longrightarrow B A$ induces a natural morphism of
supercomplexes
$$\hat{X}(C)\otimes \hat{X}(D)\longrightarrow \hat{X}(B A),\tag 9$$
where $X$ is the $X$ complex functor of Cuntz and Quillen. We will then
apply this result to the structure map of a homotopy $G$ algebra to
construct operations on the cyclic and Hochschild homology of homotopy 
$G$ algebras.

Note that, in general,
there is no natural  map  $X(C)\otimes X(D) \longrightarrow
X(C\otimes D)$; otherwise defining (9) would be a trivial matter.
This is simply because $X$ only captures homological information up to 
dimension one.
Instead, we obtain (9) as a composition
$$\hat{X}(C)\otimes \hat{X}(D)\longrightarrow \hat{X}^2(C\otimes D)
\longrightarrow \hat{X}^2(B
A)\longrightarrow \hat{X}(B A), $$
where $X^2$ is a certain refinement of $X$ to capture degree 2 homology 
classes. Despite the fact that there is no natural transformation 
$X^2\longrightarrow X$, we can however make use of the fact that the 
underlying coalgebra of $BA$ is free and show that $\hat{X}(BA)$ is a
deformation retract of $\hat{X}^2(BA)$. This gives the last map  in the
above sequence. The
first map is simply the $DG$ coalgebra analogue of a map constructed by 
M. Puschnigg in his study of Kunneth formulas in cyclic homology [13].

We need to adopt some basic definitions and constructions from [9] to
our $DG$ coalgebraic set up. Let $C$ be a $DG$ coalgebra  and let $(\Omega C,d)$ denote the $DG$
coalgebra  of {\it universal codifferential forms } over $C$. Let $\eta
:C\longrightarrow k$ be the counit of $C$. We have $\Omega^n C= C\otimes
\bar{C}^{\otimes n}$, where $\bar{C}=Ker \eta$. Let
$b:\Omega^{\bullet}C\longrightarrow \Omega^{\bullet +1}C$
be the analogue of the {\it Hochschild boundary operator} and let
$N$ be the {\it number operator} which multiplies a differential form by 
its degree. Let
$$\Omega^{norm} C=ker\{(b+dN)^2 :\Omega C\longrightarrow \Omega C\}.$$

Equipped with the differential $b+dN$ and with its natural ${\Bbb Z}/2$
grading, $(\Omega^{norm} C,b+dN)$ can be regarded as a supercomplex. There is a
decreasing filteration $\{F^n \Omega^{norm} C\}_{n\geq 2}$ on 
$\Omega^{norm} C$, where $F^n$ consists of forms of degree at least $n$. 
The successive quotient complexes $\Omega^{norm} C/F^n$ approximate the
normalized cyclic bicomplex for $DG$ coalgebras. We need only the first two quotients,
denoted by $X(C)$ and $X^2(C)$. These are the supercomplexes
$$
\align
X(C)&: \quad \quad \quad \quad \quad \quad \quad
C\overset {\overset b\to \longrightarrow } \to {\underset d\to
\longleftarrow}
 \Omega^1 C_{\natural}\\
X^2(C)&:\quad \quad \quad \quad \quad \quad C\bigoplus \Omega^2C_{\natural}
\overset {\overset b+2d  \to
\longrightarrow}\to {\underset b+d \to \longleftarrow }
\dot{ \Omega}^1 C,
\endalign
$$
where $\natural$ denotes the cocommutator subspace and $\dot{\Omega}^1 C
=\Omega^{norm ,1} C$. Note that $\Omega^1 C_{\natural}\subset
\dot{\Omega}^1 
C$.

We use $\partial_1$ (resp. $\partial_2$) to denote the horizontal (resp.
vertical) differentials in $X(C)$ and $X^2(C)$. We are mostly interested in the
total complexes of these bicomplexes which we denote by $\hat{X}(C)$ and
$\hat{X}^2(C)$. We express an even, or odd, element of $\hat{X}(C)$ as
$\omega_0+\omega_1$, where $\omega_0 \in C$ and $\omega_1\in
\Omega^1C_{\natural}$.  Similarly we write $\omega_0+\omega_2+\omega_1$
to denote an even, or odd, element of $\hat{X}^2 (C)$.
Note that we have a natural morphism of supercomplexes
$$I:\hat{X}(C)\longrightarrow \hat{X}^2(C),$$
obtained from the inclusions $C\longrightarrow C\oplus \Omega^2 
C_{\natural}$ and $\Omega^1C_{\natural}\longrightarrow \dot{\Omega}^1 C.$

In general, there is no natural map $\hat{X}^2(C)\longrightarrow 
\hat{X}(C)$. However, we would like to show that if $C=BA$ is the bar
construction, then $I$
is a homotopy equivalence and find an explicit homotopy inverse
$R:\hat{X}^2(C)\longrightarrow \hat{X}(C)$.  The easiest way to find $R$ 
is to apply homological perturbation theory. Indeed a simple version of
the so called
 perturbation lemma which we recall now is enough for our purpose.

Recall that a (super)complex $(L,\partial_1)$ is a {\it special deformation
retract} of a (super)complex $(M,\partial_1)$ if there are chain maps
 $$L\overset i\to \longrightarrow M\overset r \to \longrightarrow L$$
  and a homotopy $h:M\longrightarrow M$, of odd degree, such that 
  $ri=1_L$, $ir=1_M+\partial_1 h+h\partial_1$ and $hi=0$. In particular
  $i$ is a homotopy equivalence and $r$ is a homotopy inverse to $i$. Let us perturb
  the differentials to $\partial_1+\partial_2$ and assume that 
  $\partial_2 i=i\partial_2$. It is natural to ask if 
  $(L,\partial_1+\partial_2)$ remains a deformation retract of $(M, 
  \partial_1+\partial_2)$. It is not difficult to show that this is indeed the case, provided the
  operator $K= \sum_{n \geq 0}  (\partial_2h)^n$ can be rigorously defined. In
  this case one shows that the chain maps
  $$L\overset I\to \longrightarrow
  M \overset R\to \longrightarrow L$$
  and the homotopy $H:M\longrightarrow M$
  defined by $I=i$, $R=rK$ and $H=hK$, provide a  special deformation
  retract of $(M, \partial_1+\partial_2)$ to $(L, \partial_1 
  +\partial_2)$. In our applications $K$ will be a finite sum.

Let $C=B A $ be the bar construction of an algebra $A$ with its  counit
$\eta :BA\longrightarrow k$ and let
$\bar{B}A= ker \eta =\bigoplus_{n\geq 1} A^{\otimes n}$.
We have $\Omega^1 B A =B A\otimes \bar{B} A\simeq
B A \otimes A\otimes B A $ , $\Omega^1 BA_{\natural}\simeq =A\otimes BA $, and
$\Omega^2 B A =B A \otimes \bar{B} A \otimes \bar{B} A $. We fix a left
inverse $\theta :\Omega^1 BA \longrightarrow \Omega^1 BA_{\natural}$
for the inclusion $\Omega^1BA_{\natural}\hookrightarrow \Omega^1BA$,
defined by
$$ \theta(\alpha \otimes a \otimes \beta )= \eta (\alpha ) a \otimes
\beta .$$

>From [9], one knows that one can use {\it connections} to construct homotopy
operators  with good algebraic properties  for the Hochschild and cyclic complex
of  $DG$ coalgebras.
 Let us define
an operator $\nabla :\Omega^2 B A \longrightarrow \Omega^1 B A $, which is
supported on $B A \otimes \bar{B} A \otimes A$, by the formula
   $$ \nabla (\beta \otimes \alpha \otimes a)=\alpha \otimes a \otimes \beta ,$$
   where $\beta \in BA ,\alpha \in \bar{B}A$, and $a\in A$.  Define an
   odd
     operator $h':\hat{X}^2(B A) \longrightarrow \hat{X}^2 (B
   A)$ via the formula
   $$h'(\omega_0 +\omega _2 +\omega_1)= \nabla \omega_2.$$
    Also define even operators $ r':\hat{X}^2 (BA)\longrightarrow
   \hat{X}(BA)$ and $i':\hat{X}(BA)\longrightarrow \hat{X}^2(BA) $ via the formulas
   $$r'(\omega_0 +\omega_2+\omega_1)=\omega_0 +\theta \omega_1 ,$$
   and $i'=I$.

\proclaim{Lemma 3} $(r',h',i')$ is a special deformation
retract of $(\hat{X}^2 (B A),b)$ to
$(\hat{X}( B A),b)$.
\endproclaim
\demo{Proof} The relations $r'i'=1$ and $h'i'=0$ are easy to verify. The
relation $ i'r'=1+bh'+h'b$ amounts to $ \theta \omega_1
=\omega_2+\omega_1 +b\nabla \omega_2 +\nabla b \omega_1$. This is
equivalent to showing that, for all $\omega_1 \in \dot{\Omega}^1 BA$ and
$\omega_2 \in \Omega^2 BA_{\natural} $,
$$
\align
\omega_1 +\nabla b \omega_1 &=\theta \omega_1\\
\omega_2 +b\nabla \omega_2 &= 0
\endalign
$$
While it is possible to prove these relations by a direct computation,
it is perhaps more instructive to prove the corresponding dual
statements for the tensor coalgebra $TA$. In this case the connection
$\nabla :\Omega^1 TA \longrightarrow \Omega^2 TA$ is given by $\nabla
(\alpha \otimes a \otimes \beta )=\beta d\alpha da$. To prove the first
relation let $\omega_1 = \alpha \otimes a \otimes \beta $. We
have
$$
\align
b\nabla \omega_1 &= b (\beta d\alpha da)=-b(da\beta d\alpha)\\
     &=[da\beta ,d\alpha]=da\beta \alpha -\alpha da \beta .
\endalign
$$
And hence
$$b\nabla \omega_1 +\omega_1 =da \beta \alpha =\theta \omega_1.$$

To prove the second relation, let $\omega_2=a_0da_1da_2$. We have
$$
\align
\nabla b\omega_2 &= -\nabla [a_0da_1 ,a_2]=-\nabla (a_0da_1 a_2 -a_2a_0
da_1)\\
&= -a_2da_0 da_1 +d(a_2a_0)da_1=da_2 a_0da_1=-a_0da_1da_2 \\
&= -\omega_2.
\endalign
$$
The lemma is proved.
\enddemo

To pass from the $b$-complex to $\partial_1$- complex, we compute the
operator
$k:\hat{X}^2 \longrightarrow \hat{X}^2$. It is given by
$$k(\omega_0+\omega_2+\omega_1)=(\omega_0+d\nabla \omega_2)+\omega_2
+\omega_1.$$
Invoking the perturbation lemma, let us now define the operators $h,r $ and
$i$ by the formulas:
$$
\align
h(\omega_0 +\omega_2+\omega_1)&= \nabla \omega_2 \\
r(\omega_0+\omega_2+\omega_1)&=(\omega_0 +d \nabla \omega_2) +\theta
   \omega_1,
\endalign
$$
and $i=I$.

\proclaim{Lemma 4} $(r,h,i)$ is a special deformation retract of
$(\hat{X}^2 (BA),
\partial_1)$ to $(\hat{X} (BA) ,\partial_1)$.
\endproclaim

We  use the perturbation lemma once again to pass from the
$\partial_1$-complex to $\partial_1+\partial_2$-complex. The operator
$K$ is now given by
$$K(\omega_0+\omega_2 +\omega_1)=\omega_0+\omega_2 +(\omega_1+\partial_2
\nabla \omega_2).$$
Let us define
 the
operators $R$ and $H$ by
$$
\align
R(\omega_0+\omega_2+\omega_1) &=(\omega_0+d\nabla \omega_2)+\theta
    (\omega_1 +\partial_2 \nabla \omega_2) \tag 10\\
H(\omega_0 +\omega_2 +\omega_1)&= \nabla \omega_2.
\endalign
$$

\proclaim{Proposition 5} $(R,H,I)$ is a special deformation retract of
$\hat{X}^2 (BA)$ to $\hat{X} (BA)$.
\endproclaim

Although we won't need it in this paper, we note that the above
proposition and its proof remain valid in the more general case where
$A$ is a $DG$ algebra.

In his study of Kunneth formulas in cyclic homology [13], M. Puschnigg
constructed a natural map $X^2(A\otimes B)\longrightarrow X(A)\otimes 
X(B)$, where $A$ and $B$ are algebras and the tensor product of 
supercomplexes is understood in the right hand side. This map lifts
Connes' external product in cyclic homology [1] to the level of chains in
the $X$ complex. It is given by
$$
\align
a_0b_0&\mapsto  a_0\otimes b_0\\
a_0b_0d(a_1b_1)&\mapsto  \frac{1}{2}a_0da_1\otimes
[b_0,b_1]_+ +\frac{1}{2}[a_0,a_1]_+ \otimes b_0db_1 \\
a_0b_0d(a_1b_1)d(a_2b_2)&\mapsto \frac{1}{2} a_0da_1a_2\otimes
b_0b_1db_2- \frac{1}{2}a_0a_1da_2\otimes b_0db_1b_2 \\
&=\frac{1}{2} a_0d(a_1a_2) \otimes b_0b_1 db_2 -\frac{1}{2}a_0a_1 da_2
\otimes b_0d(b_1b_2),
\endalign
$$
where $[a,b]_+ =ab+ba $ and, to keep the notation simple, we have supressed the tensor product
sign on the left hand side.
This map is  functorial and can be dualized to a $DG$ coalgebra
context to define a morphism of supercomplexes
$$p: \hat{X}(C)\otimes
\hat{X}(D)\longrightarrow \hat{X}^2 (C\otimes D),$$
 where $C$ and $D$ are
$DG$ coalgebras. Let $p_{i,j}=p\vert X_i(C)\otimes X_j(D)$. We have
$$
\align
p_{0,0}&= id\\
p_{0,1}&=\frac{1}{2} R_{2,3}\circ (1\otimes (\Delta +R_{1,2}\Delta))   \\
p_{1,0}&=\frac{1}{2}R_{2,3}\circ ((\Delta +R_{1,2}\Delta )\otimes 1)   \\
p_{1,1}&=\frac{1}{2}R \circ (1\otimes \Delta \otimes \Delta \otimes
1 -\Delta \otimes 1 \otimes 1
\otimes \Delta ).
\endalign
$$
In the above formulas $R_{i,j}$ denotes the signed exchange between $i$ and $j$
factors in a tensor product and $R=R_{4,5}R_{3,5}R_{2,3}$.  This
observation, coupled with the above proposition,  proves the following
\proclaim{Theorem 6} Let $C$ and $D$ be $DG$ coalgebras and let $A$ be an
algebra. Then any morphism of $DG$ coalgebras
$C\otimes D
\longrightarrow B A$ induces a mprphism of supercomplexes
$$\hat{X}(C)\otimes \hat{X}(D)\longrightarrow \hat{X}(B A).$$
\endproclaim

In the remainder of this section we will apply this result to homotopy 
$G$ algebras. So let
 $V= \bigoplus_{i\geq 0} V_i $ be a homotopy $G$ algebra with
 structure map $\cup :BV\otimes BV \longrightarrow BV$.
Using the inclusion $V_0\longrightarrow V$ and the surjection
$V\longrightarrow V_0$, we obtain a morphism of coalgebras
$\cup_1 :BV\otimes BV_0 \longrightarrow BV_0$.
It is given by
$$
(x_1,\cdots ,x_m)\cup_1 (a_1,\cdots ,a_n)=\underset{i_1\leq i_2\leq
\cdots \leq i_m\leq n } \to \sum (-1)^{\epsilon}
\tag 11$$
$$
  (a_1,\cdots ,a_{i_1} ,x_1\{a_{i_1+1}, \cdots ,a_{i_1+d_1}\},\cdots
   ,x_m \{a_{i_m+1},\cdots ,a_{i_m+d_m}\},\cdots ,a_n),
 $$
where, $\epsilon =\sum_{p=1}^{m}(\vert x_p\vert -1 )i_p$ and $d_i=\vert x_i \vert $. In the special case
where $V=C^{\bullet}(A,A)$, the braces in (11) are given by
 $x_{i_k}\{a_{i_k+1}, \cdots ,a_{i_k+d_k}\}
=x_{i_k}(a_{i_k +1},\cdots a_{i_k+d_k})$.

\proclaim{Lemma 7} $\cup_1$ is a morphism of $DG$ coalgebras.
\endproclaim
\demo{Proof} Note that the surjection $\pi :BV\longrightarrow BV_0$ is a 
$DG$ coalgebra map. For $\alpha \in BV$ and $\beta \in BV_0$, we have
 $$\align
 b'(\alpha \cup_1 \beta )&= b'\pi (\alpha \cup \beta
 )=\pi (b'+\delta)(\alpha \cup \beta )\\
  &=\pi [(b'+\delta)\alpha \cup \beta  +(-1)^{\vert \alpha \vert}\alpha \cup
  (b'+\delta)\beta]\\
  &= (b'+\delta)\alpha \cup_1 \beta +(-1)^{\vert \alpha \vert}[\pi (\alpha \cup b'\beta )
  +\pi (\alpha
  \cup \delta \beta )]\\
  &=  (b'+\delta)\alpha \cup_1 \beta +(-1)^{\vert \alpha \vert}\alpha
  \cup_1 b'\beta ,\\
 \endalign
 $$
since $\pi (\alpha \cup \delta \beta )=0$ as $\delta \beta $ has degree one. The
 lemma is proved.
\enddemo

Let
$$P: \hat{X}(BV)\otimes \hat{X}(BV_0) \longrightarrow \hat{X}^2 (BV_0)$$
denote the composition $ P=\cup_1 p$.
Now we can apply theorem 6 to obtain
\proclaim{Theorem 8}
Let $V$ be a homotopy $G$ algebra. Then there are natural maps of
supercomplexes
$$ \hat{X}(BV)\otimes \hat{X}(BV) \longrightarrow \hat{X}(BV),$$
$$ \hat{X}(BV)\otimes \hat{X}(BV_0)\longrightarrow \hat{X}(BV_0).$$
\endproclaim

We beleive, although have not checked it, that the first product is
homotopy associative and in fact there should exist a full structure of
higher homotopis in the sense of $A_{\infty}$-algebras. The same should
be true on the corresponding pairing between homologies.

\subhead 4 \quad Higher Operations on Cyclic Bicomplex \endsubhead

Our goal in this sestion is to find explicit formulas for the second map
in theorem 8
 and to relate it to the Hochschild and cyclic homology of the homotopy $G$
algebra $V$. In the last part we apply our formulas to a very special
homotopy  $G$ algebra, namely the deformation complex of an algebra, to
obtain higher homotopy formulas in the cyclic and $b,B$ complex of the
algebra. The computations in this section are based on results and ideas
from [9].

In [12], Nest and Tsygan have defined two different types of operations
on cyclic homology. On one hand, they have defined an action of the $b,
B$ complex of the deformation complex (as a $DG$ algebra) of an algebra
$A$ on the $b,B$ complex of $A$. This corresponds to the pairing (14)
below,
though we do not know to what extent the explicit formulas match. As
is explained in [12], this operation is very general and in particular
yields Cartan homotopy formulas for the action of higher Hochschild
cochains on the $b,B$ complex. Secondly, they have defined an action of
the Chevalley-Eilenberg complex of the deformation complex (as a $DG$
Lie algebra) on the $b, B$ complex of $A$. Most probably, this operation
too is a consequence of the homotopy $G$ algebra structure of the
deformation complex by a similar pattern as we derived (14) from  theorem 8.

Let $\Omega_0 =(x_1,\cdots ,x_m) \in BV$, $\Omega_1 =y_0\otimes
(y_1,\cdots ,y_n) \in \Omega^1 BV_{\natural} $, $\omega_0=(a_1,\cdots
a_p) \in BV_0 $ and $\omega_1 =b_0\otimes (b_1,\cdots b_q) \in \Omega^1
BV_{0 \natural}$. Let $\eta_0 +\eta_1 \in \hat{X}(BV_0)$ be the image of
$(\Omega_0 +\Omega_1)\otimes (\omega_0 +\omega_1)$ under the second map
in
Theorem 8. Using  formula (10) for the retraction $R$, we have
$$
\align
\eta_0 &= P(\Omega_0 \otimes \omega_0) +d\nabla
P(\Omega_1 \otimes \omega_1)\\
\eta_1 &= \theta P(\Omega_0 \otimes \omega_1 +\Omega_1 \otimes
\omega_0)+\theta \partial_2 \nabla  P(\Omega_1\otimes \omega_1).
\endalign
$$
Note that
$$  P(\Omega_0\otimes \omega_0)= (x_1,\cdots ,x_m)\cup_1
(a_1,\cdots ,a_p).$$

To simplify the notation, we resort to the following convention.
It is better to consider the points $y_0, y_1, \cdots , y_n$ as located
in the clockwise order on the circle. Let $\pi_k$
denote the set of all partitions of these points on the circle
into $k$ intervals. We allow one or several of these interavls to be
empty, in which case they represent 1. For
 example, for $n=2$,   $\pi_1 $ has 3 elements while $\pi_2$ has 12
 elements. We denote an element of $\pi_2$ by a pair $(\alpha ,\beta )$,
 and similarly for elements of $\pi_k$. It is also convenient
 to write $X(Y)$ for $ X\cup_1 Y $.

 To compute the other components of $\eta_0$ and $\eta_1$,
first we have to find the image of $\Omega_1$ under the inclusion
$\Omega^1 BV_{\natural} \longrightarrow \Omega^1 BV$. We have
$$
\align
y_0\otimes (y_1,\cdots ,y_n)&\mapsto
\sum_{i=0,j=0}^{n,n-i}(-1)^{\epsilon}(y_{i+1},\cdots ,y_{i+j})\otimes (y_{i+j+1},\cdots
y_{i})\\
 &=\sum_{\underset{y_0\in \beta} \to {(\alpha , \beta )}}(-1)^{\vert
 \alpha \vert \vert \beta \vert }\alpha \otimes
 \beta
 \endalign
 $$
and similarly for $ \omega_1$. Now we  have
$$
\align
d\nabla  P(\Omega_1\otimes \omega_1) &= \frac{1}{2} \sum_{i=0}^{n}
\sum_{\underset{y_0\in (\beta ,y_i)}\to {(\alpha ,\beta ,y_i)}}
\sum_{\underset{b_0\in \gamma'}\to {(\alpha' , \beta' ,\gamma')
}}  \pm (\beta (\beta'), y_i(\gamma') ,\alpha (\alpha'))
\tag 12 \\
&   -\frac{1}{2} \sum_{(\alpha ,\beta ,y_0)}\sum_{\underset{b_0\in
(\beta',\gamma')}\to {(\alpha',\beta', \gamma')}}\pm (\beta
(\beta'), y_0(\gamma'),\alpha (\alpha')).
\endalign
$$
The signs can be easily made explicit.
Next we compute the contribution of $\Omega_0 \otimes \omega_1
+\Omega_1 \otimes \omega_0 $ to $\eta_1$. First note that
$\theta (\alpha \otimes \beta) \neq 0$ only if $\alpha =1 $ or $\alpha
\in V_0$. In this case we  have
$$
\align
\theta (a_0\otimes (a_1,\cdots a_n)) &=a_0\otimes (a_1,\cdots a_n),\\
\theta (1\otimes (a_1,\cdots ,a_n))&=-a_1\otimes (a_2,\cdots ,a_n).
\endalign
$$
Using the above information  plus our formulas for P, we obtain
$$
\align
\theta  P(\Omega_0 \otimes \omega_1)&=- \sum_{(\alpha')} (x_1,\cdots ,x_m)
(\alpha')\\
&+\frac{1}{2} \sum_{\underset{b_0 \in \beta' } \to { (\alpha' ,\beta')}}
\pm (x_1 (\alpha' ), (x_2,\cdots ,x_m) (\beta'))\\
&+\frac{1}{2} \sum_{\underset{b_0 \in \beta'}\to {(\alpha' ,\beta')}}
\pm (x_m
(\alpha') ,(x_1,\cdots ,x_{m-1})(\beta')).
\endalign
$$
Similarly we have
$$
\align
\theta  P(\Omega_1 \otimes \omega_0)&= \frac{1}{2} \sum_{i=1}^{n}
\sum_{(y_i,\beta )} \pm (y_i(a_1,\cdots ,a_{d_i}), \beta (a_{d_i +1},\cdots
,a_p))\\
&- \sum_{(\beta ) }  \beta (a_1,\cdots , a_p )\\
&+ \frac{1}{2}\sum_{i=1}^n \sum_{(y_i,\beta )} \pm (y_i(a_{p-d_i}, \cdots ,
a_p),\beta (a_1,\cdots ,a_{p-d_i-1})),
\endalign
$$
  where $d_i=\vert y_i\vert $.

Finally let us  compute the contribution of $\Omega_1 \otimes \omega_1$ to
$\eta_1$. Using the formulas for the connection $\nabla$ and for the induced
differential $\partial_2$  we note that $\theta \partial_2 \nabla (\beta
\otimes \alpha \otimes \gamma )\neq 0$ only if $\alpha \in V_0$ and
$\gamma \in V_0$. In this case we have
   $$\theta \partial_2 \nabla (\beta \otimes \alpha \otimes \gamma
   )=\alpha \gamma  \otimes
   \beta .$$
>From this we obtain
$$
\align
\theta \partial_2 \nabla P(\Omega_1\otimes \omega_1)&=
 \frac{1}{2} \sum_{\underset{b_0\in \gamma'}\to  {(\alpha' ,\beta'
 ,\gamma' )}}\pm (y_0 (\beta') y_1
 (\gamma') , (y_2, \cdots ,y_n) (\alpha') ) \tag 13 \\
 &-
 \frac{1}{2} \sum_{\underset{b_0\in \beta'}\to {(\alpha' , \beta',
 \gamma') }} \pm (y_n (\beta') y_0
 (\gamma') , (y_1,\cdots
y_{n-1}) (\alpha')).
\endalign
$$

We define the {\it periodic cyclic homology}
of a $DG$ algebra $V= \bigoplus_{i\geq 0} V_i$ to
be the homology of the supercomplex $\hat{X}(BV)$. In the special case
when $V=A$ is an algebra, the bicomplex $X(B A)$ has been shown by
Quillen [14] to be isomorphic, up to a shift in the vertical direction,
to the {\it cyclic bicomplex} of $A$. The same
argument works in the $DG$ case. Let $\Cal C (V)$ denote the total
complex of the cyclic bicomplex of $V$. We have
$$\Cal C_{ev}(V) =\Cal C_{odd}(V)=\prod_{n\geq 0}C_n(V),$$
 where
 $C_n(V)=V^{\otimes (n+1)}$.
 Theorem 8 can now be reformulated as
 \proclaim{Theorem 9}
  Let $V$ be a homotopy $G$ algebra. Then there is a natural
morphism of supercomplexes
$$\Cal C (V) \otimes \Cal C(V_0)\longrightarrow \Cal C (V_0).$$
\endproclaim

Finally we turn to Connes' $b, B$  bicomplex  and the analogue of our
formulas in that context. This is important because in many applications
of cyclic homology and noncommutative geometry [1, 2] the $b,B$ bicomplex appears
in a natural way. Using an
explicit homotopy equivalence between the cyclic and $b,B$ bicomplexes,
we can transform theorem 9 into a morphism between $b,B$ complexes.
One obtains, however, simpler formulas if
one restricts to the {\it normalized} $b, B$ complex.

Let $(V,\delta )$ be a unital $DG$ algebra and let $\Cal B (V)$ denote its $b,B$
complex. We have
$$\Cal B (V)_{ev}=\prod_{n\geq 0} C_{2n} (V) \quad \quad \quad  \quad
\Cal B (V)_{odd}=\prod_{n\geq 0}C_{2n+1}(V)$$
The differential is
given by $b+B+\delta $, where $B:C_{\bullet}(V)\longrightarrow
C_{\bullet +1}(V)$ is Connes' boundary operator. Let $\bar{V}=V/k$ and
$\bar{C}_n(V)=V\otimes \bar{V}^{\otimes n}$. The normalized $b,B$
complex of $V$, denoted $\bar{\Cal B}(V)$,  is defined similarly except that we replace $C_n(V)$ by
$\bar{C}_n(V)$.

For $n\geq 0$, let $\lambda:C_n(V) \longrightarrow C_n(V)$ denote the cyclic
shift operator, let $N$ be the correspnding
norm operator and let $s:C_n (V)\longrightarrow C_{n+1}(V)$ be defined
by
$s(v_0,\cdots ,v_n)=(1,v_0,\cdots ,v_n)$.
Recall the morphisms of complexes
$$I: \Cal B (V) \longrightarrow \Cal C (V) \quad \quad \quad \quad
J: \Cal C (V) \longrightarrow \Cal B (V)$$
defined  by
$$
I= 1+sN   \quad \quad \quad \quad \quad \quad
 J=1+s(1-\lambda)
 $$
It is known that the
operators $I$ and $J$ are homotopy
inverse to each other [10].

Using the chain maps $I$ and $J$ it is clear that we can transform
theorem 9 into a morphism of
complexes
$$ \Cal B (V)\otimes \Cal B (V_0)\longrightarrow \Cal B (V_0).$$

We specialize to the case where $V=C(A,A)$ is the deformation complex of
a unital algebra $A$.  Note that in this case
$V_0=C^0(A,A)=A$. A  cochain $\phi \in C^n(A,A)$ is said to be {\it normalized} if $\phi (a_1,\cdots
,a_n)=0$ whenever $a_i=1$ for some $i$. Let $C_{norm}(A,A)$ denote the
subcomplex of normalized cochains. It is easy to check that
$C_{norm}(A,A)$ is indeed a sub $DG$ algebra of the $DG$ algebra
$C(A,A)$. Hence we have an inclusion
$ \Cal B (C_{norm}(A,A))\longrightarrow  \Cal B (C(A,A))$, and a
morphism of supercomplexes
   $$\Cal B (C_{norm}(A,A))\otimes \Cal B (A)\longrightarrow \Cal B
   (A).$$
Now our explicit formulas show that the above map descends to define
a morphism of supercomplexes
$$\Cal B (C_{norm}(A,A)) \otimes \bar{\Cal B}(A)\longrightarrow
\bar{\Cal B}(A).\tag 14$$

We denote this map as well as the one in theorem 9  by $\cup$. We
obtain explicit formulas for $\cup$ as follows. Let $D=(D_0, \cdots ,D_m) \in
{\Cal B}(C_{norm} (A,A))$ and $ a =(a_0,\cdots , a_n) \in \bar{{\Cal B}}(A).$ Let
$ID=\Omega_0 +\Omega_1 =sND+D$ and $I a =\omega_0 +\omega_1
=sNa + a .$ Also let $ID\cup I a =\eta_0 +\eta_1 $. We have
$$D\cup a = J(ID\cup I a )=J(\eta_0 +\eta_1)=
s(1-\lambda)\eta_0+\eta_1.$$
 Now we have
 $$s(1-\lambda)\eta_0=s(1-\lambda)(P(sND \otimes sN a )+d\nabla
 P(D\otimes a )).$$
 Because of our normalization conditions  we have
       $$s(1-\lambda)P(sND \otimes sN a )=0 .$$
 Using (12) we get
 $$
 \align
 s(1-\lambda )\eta_0 &= s(1-\lambda )d\nabla P(D\otimes a)\\
                     &= \frac{1}{2}\sum_{i=0}^{m}\sum_{\underset{D_0\in
                     (\beta , D_i)}\to {(\alpha ,\beta
                     ,D_i)}}\sum_{\underset{a_0 \in \gamma'}\to {(\alpha',
                     \beta' , \gamma')}} \pm (\beta (\beta'), D_i (\gamma')
                     ,\alpha (\alpha')))\\
                     & -\frac{1}{2}\sum_{(\alpha ,\beta
                     ,D_0)}\sum_{\underset{D_0\in (\beta',
                     \gamma')}\to{(\alpha', \beta' , \gamma')}}\pm
                     (\beta (\beta'), D_0(\gamma'),\alpha (\alpha')).
 \endalign
$$

Similarly  we have
$$
\eta_1= \theta P(sND \otimes \alpha + D\otimes sN a )+\theta
\partial_2 \nabla P(D\otimes a ).$$

\proclaim{Lemma 10}
We have
 $\theta P (D\otimes sN a )=0$ and
$\theta P(sND \otimes a )= 0.$
 \endproclaim

Using the above lemma and (13), we get
$$
\align
\eta_1 &= \theta \partial_2 \nabla P( D\otimes a ) \\
       &= \frac{1}{2}\sum_{\underset{a_0 \in \gamma'}\to {(\alpha'
       ,\beta', \gamma')}}\pm (D_0(\beta')D_1(\gamma'), (D_2,\cdots
       ,D_n)(\alpha'))\\
       &- \frac{1}{2}\sum_{\underset{a_0\in \beta'}\to {(\alpha',
       \beta', \gamma')}}\pm (D_n(\beta')D_0(\gamma'), (D_1,\cdots
       ,D_{n-1})(\alpha')).
       \endalign
       $$

Note that because of homogeneity we can drop the unpleasant factor of
$\frac{1}{2}$ from the formulas. In the unnormalized case, however, this
can not be done.

\ref\key 15
\by M. Kontsevich
\paper Deformation quantization of Poisson manifolds
\jour q-alg/9709040
\yr 1997
\endref

\ref\key 16
\by M. Kontsevich
\paper Private communication
\endref
\enddocument